\documentclass[12pt]{article}
\usepackage[expansion=false]{microtype}
\usepackage[reqno]{amsmath}
\usepackage{fixmath}
\usepackage{amssymb, amsthm, enumerate}
\usepackage{mathtools}
\usepackage{graphicx}

\usepackage{hyperref}

\DeclarePairedDelimiter{\abs}{\lvert}{\rvert}
\DeclarePairedDelimiter{\norm}{\lVert}{\rVert}

\newtheoremstyle{plainsl}%
    {\topsep}
    {\topsep}
    {\slshape} 
    {}
    {\normalfont\bfseries}
    {.}
    { }
    {}

\swapnumbers

{\theoremstyle{plainsl}
\newtheorem{theorem}{Theorem}[section]
\newtheorem{lemma}[theorem]{Lemma}
\newtheorem{corollary}[theorem]{Corollary}}
{\theoremstyle{remark}
}

\renewcommand\proof{\noindent\textsl{Proof. }}
\newcommand\sqr[2]{{\vbox{\hrule height.#2pt
    \hbox{\vrule width.#2pt height#1pt \kern#1pt
        \vrule width.#2pt}\hrule height.#2pt}}}
\renewcommand\qed{%
    \ifmmode\eqno\sqr53
    \else\nolinebreak\ \hfill\sqr53\medbreak\fi}

\numberwithin{equation}{section}

\newcommand\de{\delta}

\newcommand\eps{\epsilon}
\newcommand\ga{\gamma}

\renewcommand\th{\theta} 

\newcommand\cA{{\mathcal A}}

\newcommand\cx{{\mathbb C}}

\newcommand\re{{\mathbb R}}
\newcommand\rats{{\mathbb Q}}

\newcommand\diff{\mathbin{\mkern-1.5mu\setminus\mkern-1.5mu}}
\newcommand\sbs{\subseteq}

\newcommand\seq[3]{#1_{#2},\ldots,#1_{#3}}

\DeclareMathOperator{\aut}{Aut}

\newcommand\pmat[1]{\begin{pmatrix} #1 \end{pmatrix}}

\DeclareMathOperator{\rk}{rk}
\DeclareMathOperator{\tr}{tr}
\DeclareMathOperator{\col}{col}
\newcommand\ip[2]{\langle#1,#2\rangle}

\newcommand\alg[1]{\langle #1\rangle}
\DeclareMathOperator\mult{mult}

\newcommand\one{\textbf{1}}

\newcommand\hmx{\widehat{M}_X}
\DeclareMathOperator{\comm}{comm}

\title{Strongly Cospectral Vertices}

\author{Chris Godsil, Jamie Smith}

\begin{document}
\maketitle

\begin{abstract}
	Two vertices $a$ and $b$ in a graph $X$ are \textsl{cospectral} if the vertex-deleted
	subgraphs $X\diff a$ and $X\diff b$ have the same characteristic polynomial. In this
	paper we investigate a strengthening of this relation on vertices, that arises in
	investigations of continuous quantum walks. Suppose the vectors $e_a$ for $a$ in $V(X)$
	are the standard basis for $\re^{V(X)}$. We say that $a$ and $b$ are 
	\textsl{strongly cospectral} if for each eigenspace $U$ of $A(X)$, the orthogonal
	projections of $e_a$ and $e_b$ are either equal or differ only in sign.
	We develop the basic theory of this concept and provide constructions of graphs
	with pairs of strongly cospectral vertices. Given a continuous quantum walk on
	on a graph, each vertex determines a curve in complex projective space. We derive
	results that show tht the closer these curves are, the more ``similar'' the 
	corresponding vertices are.
\end{abstract}

\tableofcontents

\section{Introduction}

To start, we set up some machinery for working with quantum states.
We will represent a quantum state in $\cx^n$ by a \textsl{density matrix}, a positive
semidefinite $n\times n$ matrix with trace one. A density matrix $D$ represents a
\textsl{pure state} if $\rk(D)=1$, in which case $D=zz^*$ for some unit vector $z$.
We will only be concerned with pure states in this paper and generally they will be
associated to vertices of a graph $X$---if $a\in V(X)$, then $e_a$ denotes the standard basis
vector in $\cx^{V(X)}$ indexed by $a$ and our focus will be on pure states of the
form $D_a=e_ae_a^T$. If $D$ is a pure state then $D^2=D$ and $D$ represents orthogonal
projection onto the column space of $D$; thus $D$ corresponds to a point in complex
projective space.

If $X$ is a graph with adjacency matrix $A$, the \textsl{continuous quantum walk} on $X$
is determined by the family of unitary matrices
\[
	U(t) = \exp(itA),\qquad t\ge0.
\]
The understanding is that if, initially our system is in the satate associated with the
density matrix $D$, then at time $t$ its state is given by
\[
	U(t)DU(-t).
\]
It is easy to check that this is a density matrix, which we denote by $D(t)$, and
that $D(t)$ is pure if and only if $D$ is. It follows that, if our initial state $D$ is pure, 
a quantum walk determines a curve in projective space, namely the set of points $D(t)$.
(If our initial state were not pure, we would have a curve on a Grassmannian, but we will not go
there.)

Given distinct vertices $a$ and $b$ in $X$, one question of interest to physicists is 
whether there is a time $t$ such that $D_b$ lies on the curve containing $D_a$; equivalently
is there a time $t$ such that $U(t)D_aU(-t)=D_b$. If there is such a time, we say that
we have \textsl{perfect state transfer} from $a$ to $b$ at time $t$. If we do have
perfect state transfer at time $t$, then
\[
	\norm{D_a(t)-D_b} = 0.
\]
Since, as it happen, perfect state transfer is rare, we might decide to settle for less:
we could ask whether, given $\eps>0$, there is a time $t$ such that
\[
	\norm{D_a(t)-D_b} < \eps.
\]
If this is possible (for all positive $\eps$) we have \textsl{pretty good state transfer}
from $a$ to $b$. Pretty good state transfer occurs more often than perfect state transfer. 
For example we get perfect state transfer between the end-vertices of the path $P_n$ if and 
only if $n=2$ or $n=3$, but we have pretty good state transfer between the end-vertices of $P_n$ 
if and only if $n+1$ is a power of two, a prime, or twice a prime. For details see Banchi et al \cite{path-pgst}; more recent work on this topic appears in \cite{gc-kg-cvb, cvb-pgstchar}.)

Let $\seq\theta1m$ be the distinct eigenvalues of the adjacency matrix $A$ of the graph $X$. 
For each eigenvalue $\theta_r$ there is an idempotent matrix $E_r$
representing orthogonal projection onto the eigenspace with eigenvalue $\theta_r$.
If $f$ is a function defined on the eigenvalues of $A$, then
\[
	f(A) = \sum_r f(\th_r)E_r
\]
and, in particular
\[
	U(t) = \sum_r e^{it\th_r}E_r.
\]
Hence
\[
	D(t) = \sum_{r,s} e^{it(\th_r-\th_s)} E_rDE_s
\]
and so $D_a(t)=D_b$ if and only if
\[
	\sum_{r,s} e^{it(\th_r-\th_s)} E_rD_aE_s = D_b = \sum_{r,s} E_rD_bE_s
\]
and this holds if and only if
\[
	e^{it(\th_r-\th_s)} E_rD_aE_s = E_rD_bE_s
\]
for all $r,s$. Now all six matrices in this equality are real, whence we deduce that if
perfect state transfer occurs,
\[
	e^{it(\th_r-\th_s)} E_rD_aE_s = \pm1
\]
and, for each $r$,
\[
	E_rD_aE_r = E_rD_bE_r.
\]
(The diagonal entries in both sides are necessarily non-negative since density matrices 
are positive semidefinite and $E_r$ is symmetric, whence both sides are positive semidefinite.) 
This leads us to the conclusion that, if perfect state transfer
from $a$ to $b$ occurs, then for each $r$.
\[
	E_re_a = \pm E_re_b
\]

We define two vertices $a$ and $b$ in a graph $X$ to be \textsl{strongly cospectral}
if, for each spectral idempotent $E_r$ of $X$, we have $E_re_a=\pm E_re_b$. Our ruminations
have lead to the conclusion that, if there is perfect state transfer between vertices $a$
and $b$, then these two vertices are strongly cospectral. There is a related
and older concept, due to Schwenk \cite{cosp-trees}: vertices $a$ and $b$ in the graph $X$
are \textsl{cospectral} if the vertex-deleted subgraphs $X\diff a$ and $X\diff b$
are cospectral. We will see that `strongly cospectral' is a refinement of this concept.
(The first explicit appearance of strongly cospectral vertices is probably in \cite{dbl-stars}.)

The first part of this paper develops the theory of strongly cospectral vertices. We show
that if vertices $a$ and $b$ in $X$ are strongly cospectral, then any automorphism of $X$
that fixes $a$ must fix the vertex $b$. (So the concept has combinatorial implications.)
We provide a number of characterizations, for example: vertices $a$ and $b$
are strongly cospectral if and only if they. are cospectral and all poles of the
rational function $\phi(X\diff\{a,b\},t)/\phi(X,t)$ are simple. We use this to provide 
constructions of graphs with pairs of cospectral vertices. We show that cospectral vertices
and strongly cospectral vertices are connected by mappings that can viewed as relaxations
of automorphisms. Thus we prove that $a$ and $b$ are strongly cospectral, there is
an orthogonal matrix $Q$, a rational polynomial in $A$, such that $Q^2=I$ and $Qe_a=e_b$.

The second part of this paper considers the geometry of the orbits of the pure states
of the form $D_a$. As we noted above, there is perfect state transfer from $a$ to $b$
if and only if $D_b$ lies in the orbit of $D_a$; equivalently if and only if the
orbits of $D_a$ and $D_b$ coincide. Further we have pretty good state transfer
if and only if $D_b$ lies in the closure of the orbit of $D_a$, that is, if and
only if the closures of the two orbits are equal. We show that equality of orbits,
or of orbit closures, is unnecessary. Among other things, we prove that if the orbits are 
sufficiently close, then $a$ and $b$ must be cospectral and, if they are even closer,
then $a$ and $b$ must be cospectral.

\section{Cospectral Vertices}

\begin{figure}[htbp]
\begin{center}	
	\includegraphics[totalheight=3cm]{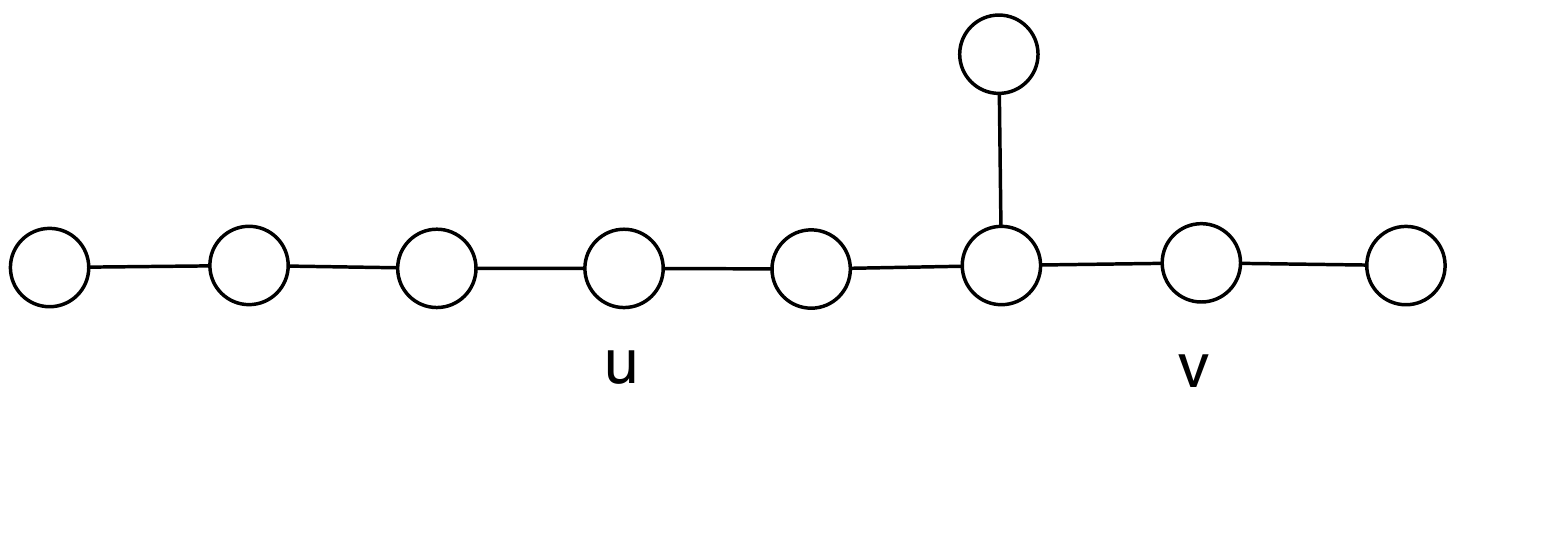}\\
	\caption{\label{fig:Fig1} A pair of cospectral vertices}
\end{center}
\end{figure}

We view the relation of being strongly cospectral as a combination of two relations.
The first of these two is an older concept: two vertices $a$ and $b$ in a graph $X$
are \textsl{cospectral} if the characteristic polynomials of the vertex-deleted subgraphs
$X\diff a$ and $X\diff b$ are equal, that is, 
\[
	\phi(X\diff a,t) = \phi(X\diff b,t).
\]
It is immediate that that if there is an automorphism of $X$ that maps $a$ to $b$,
then $a$ and $b$ are cospectral. Cospectral vertices were first introduced in Schwenk's 
fundamental paper \cite{cosp-trees}; here Schwenk noted 
that the vertices $u$ and $v$ in the tree in Figure~\ref{fig:Fig1} are cospectral, but lie in
different orbits of the automorphism group of the tree. Using this he was able to show that
the proportion of trees on $n$ vertices that are determined by their characteristic
polynomial goes to zero as $n\to\infty$.

There are a surprising number of characterizations of cospectral vertices. We will
list them in the next section, but we need first to introduce more terminology.

Suppose $S$ is a subset of the vertices a graph $X$ with characteristic vector $z$
and $n=|V(X)|$. We define the \textsl{walk matrix} $M_S$ relative to $S$ to be the $n\times n$
matrix with the vectors
\[
	z,\ Az,\ldots, A^{n-1}z
\]
as its columns. The case of interest to us will be when $S$ is a single vertex $a$ and, 
in this case, we will refer the walk matrix relative to $a$. We will use $e_S$ to denote
the characteristic vector of $S$.
The column space of $M_S$
is $A$-invariant, and so it is a module over the ring $\re[A]$ of real polynomials in $A$.
It is in fact a cyclic module, generated by the first column $z$ of $M_S$. We call it
the \textsl{walk module} relative to $S$.

We see that the $ij$-entry of of $M_S^TM_S$ is $z^TA^{i+j-2}z$, and so it is equal to the
number of walks on $X$ with length $i+j-2$ that start and end on a vertex in $S$. Hence
if $S=\{a\}$, then this entry is the number of closed walks in $X$ that start at $a$
and have length $i+j-2$. We define $W_S(X,t)$ to be the generating function
\[
	\sum_{k\ge0}z^TA^kz t^k = z^T(I-tA)^{-1}z.
\]

\begin{lemma}\label{lem:mcoeffs}
	Let $a$ and $b$ be vertices in $X$. Then $W_a(X,t)=W_b(X,t)$ if and only if
	$M_a^TM_a=M_b^TM_b$.
\end{lemma}

\proof
It should be clear that, if the walk-generating functions are equal, the matrix products are equal.
For the converse, let $\seq\th1m$ denote the distinct eigenvalues of $A$ and let E$\seq E1m$ denote 
the corresponding orthogonal projections onto the distinct eigenspaces of $A$. Then for any vector $z$,
\[
	z^T(I-tA)^{-1}z = \sum_r \frac{z^TE_r z}{1-t\th_r}.
\]
Since $m\le n$, it follows that the generating function $z^T(I-tA)^{-1}z$ is determined by its
first $m$ coefficients.\qed

\section{Characterizing Cospectral Vertices}

We give a comprehensive list of characterizations of cospectral vertices. The first four appear
already in \cite{cgbdm-wlkreg}

\begin{theorem}\label{thm:cospvs}
	Let $a$ and $b$ be vertices in the graph $X$ with corresponding walk matrices 
	$M_a$ and $M_b$. The following statements are equivalent:
	\begin{enumerate}[(a)]
		\item 
		$a$ and $b$ are cospectral.
		\item
		$\phi(X\diff a,t) = \phi(X\diff b,t)$.
		\item
		$W_a(X,t)=W_b(X,t)$.
		\item
		For each spectral idempotent $E_r$ we have $(E_r)_{a,a}=(E_r)_{b,b}$.
		\item
		For any non-negative integer $k$ we have $(A^k)_{a,a}=(A^k)_{b,b}$.
		\item
		$M_a^TM_a = M_b^TM_b$.
		\item
		The $\re[A]$-modules generated by $e_a-e_b$ and $e_a+e_b$ are orthogonal
		subspaces of $\re^{V(X)}$.
	\end{enumerate}
\end{theorem}

\proof
Claims (a) and (b) are equivalent, because (b) is the definition of cospectral.
From the proof of Lemma~\ref{lem:mcoeffs} we have
\[
	t^{-1}W_v(X,t^{-1}) = \frac{\phi(X\diff v,t)}{\phi(X,t)}
\]
and, from \cite[p.~30]{cg-blue},
\begin{equation}
	\label{eq:phiXvphiX}
	\frac{\phi(X\diff v,t)}{\phi(X,t)} = \sum_r \frac{(E_r)_{v,v}}{t-\th_r}.
\end{equation}
Hence (b), (c) and (d) are equivalent. Since any power of $A$ is a linear combination
of the spectral idempotents $E_r$, and since the spectral idempotents are polynomials in $A$,
we see that (d) and (e) are equivalent. By the discussion in the previous section, (c)
and (f) are equivalent.

We turn to (g). The given modules are orthogonal if and only if for all non-negative $i$ and $j$,
we have
\[
	\ip{A^i(e_a-e_b)}{A^j(e_a+e_b)} = 0,
\]
equivalently if and only if 
\[
	(e_a-e_b)^TA^k(e_a+e_b)=0
\]
for all $k\ge0$. This is equivalent in turn to
\[
	(e_a-e_b)^TE_t(e_a+e_b)=0
\]
for each spectral idempotent $E_r$. As
\[
	(e_a-e_b)^TE_t(e_a+e_b) = e_a^TE_re_a - e_b^TE_re_b - e_b^TE_re_a + e_a^TE_re_b
\]
and 
\[
e_b^TE_re_a = (E_r)_{b,a} = (E_r)_{a,b} =e_b^TE_re_b
\]
we find that $(e_a-e_b)^TE_t(e_a+e_b)=0$ for all $r$ if and only if $e_a^TE_re_a=e_b^TE_re_b$ 
for all $r$.\qed

We make some remarks. One consequence of Part (g) of the theorem is that if two vertices of $X$ are
cospectral, then the characteristic polynomial of $X$ factors non-trivially over $\rats$.
More precisely, the characteristic polynomials of the respective restrictions of $A$ to
the modules generated by $e_a-e_b$ and $e_a+e_b$ are disjoint factors of $\phi(X,t)$.

A graph is said to be \textsl{walk regular} if for each non-negative integer $k$,
the diagonal of $A^k$ is constant or, equivalently if the diagonals of the spectral
idempotents are constant. In a walk-regular graph, any two vertices are cospectral;
in particular any two vertices of a strongly regular graph are cospectral.

Finally, since $E_r = E_R^TE_r$, we have
\[
	(E_r)_{v,v} = e_v^TE_r^TE_re_v = \norm{E_re_v}^2,
\]
whence vertices $a$ and $b$ are cospectral if and only if the eigenspace projections
$E_re_a$ and $E_re_b$ have the same length for each $r$. It follows (as we would hope)
that strongly cospectral vertices are cospectral.

\section{Parallel Vertices: Characterizations}

We have developed some of the theory of cospectral vertices and noted that strongly
cospectral vertices are cospectral. To characterize strongly cospectral vertices,
we need a second condition. Two vertices $a$ in $b$ in $X$ are \textsl{parallel} if,
for each $r$, one of the vectors $E_re_a$ and $E_re_b$ is a scalar multiple of the other.
Equivalently $a$ and $b$ are parallel if and only if the vectors $E_re_a$ and $E_re_b$
are parallel for each $r$. As an immediate consequence of the definition of strongly
cospectral vertices, we have:

\begin{lemma}
	Two vertices in a graph are strongly cospectral if and only if they are cospectral
	and parallel.\qed
\end{lemma}

If the eigenvalues of $X$ are all simple, it is easy to see that any two vertices
in $X$ are parallel. It follows in this case that two vertices are strongly cospectral
if and only if they are cospectral. (The eigenvalues of Schwenks tree in Figure~\ref{fig:Fig1})
are simple, and the vertices $u$ and $v$ there are strongly cospectral.)

\begin{lemma}
	The eigenvalues of $X$ are all simple if and only if any two vertices of $X$
	are parallel.
\end{lemma}

\proof
Suppose any two vertices of $X$ are parallel. If $a\in V(X)$ and $E_re_a\ne0$, then for each $r$
we have that $E_re_b$ is a scalar multiple of $E_re_a$. Hence $E_re_a$ spans the eigenspace
belonging to $\th_r$ and so $\th_r$ has multiplicity one.\qed

We use $\alg{e_u}_A$ to denote the $\re[A]$-module generated by $e_u$, and we call it the 
\textsl{walk module} relative to $u$. (When $A$ is clear from the context, we may be lazy and write simply $\alg{e_u}$.) The \textsl{eigenvalue support} of a subset $S$ of $V(X)$ with 
characteristic vector $z$ is the set of eigenvalues $\th_r$ such that $E_rz\ne0$. 
(We will also refer to the 
eigenvalue support of an arbitrary vector.) Two cospectral vertices  
necessarily have the same eigenvalue support.

\begin{lemma}
	The walk modules generated by vertices $a$ and $b$ in $X$ are equal if and only if $a$ and $b$
	are parallel and have the same eigenvalue support.
\end{lemma}

\proof 
If $u\in V(X)$, the non-zero vectors $E_re_u$ form an orthogonal basis for $\alg{e_u}$. Given this,
the result is immediate.\qed

Finally we note that, by \cite[Lemma~13.1]{st-transfer}, if we have pretty good state transfer from 
vertex $a$ to vertex $b$, then $a$ are $b$ are strongly cospectral. (This result is a 
private communication from Dave Witte Morris.) Since perfect state transfer 
can be viewed as a special case of pretty good state transfer, it follows that vertices 
involved in perfect state transfer are necessarily strongly cospectral. (This is not hard
to prove directly.)

\section{Average States}

If $\seq\th1m$ are the distinct eigenvalues of the adjacency matrix $A$ of $X$, 
we use $E_r$ to denote the matrix representing orthogonal projection onto the 
$\th_r$-eigenspace of $A$. So $A$ has spectral decomposition
\[
	A = \sum_r \th_r E_r,
\]

We make use of some theory developed in \cite{us-avmix2017}. The \textsl{commutant} $\comm(A)$
of a matrix $A$ is the set of all matrices that commute with $A$. If $A$ is $n\times n$,
then $\comm(A)$ is a subspace of the space of $n\times n$ real matrices. This latter
space is an inner product space, with inner product 
\[
	\ip{M}{N} = \tr(M^TN).
\]
The $\norm{M}$ of a matrix $M$ is $\ip{M}{M}^{1/2}$. 
The operation of orthogonal projection onto $\comm(A)$ is well defined, we
denote the orthogonal projection of a matrix $M$ onto $\comm(A)$ by $\Phi(M)$.

\begin{lemma}
	If $A$ is a symmetric matrix with spectral idempotents $\seq E1m$, then
	\[
		\Phi(M) = \sum_r E_rME_r.\qed
	\]
\end{lemma}

As $\Phi$ is linear and self-adjoint,
\[
	\ip{\Phi(M)}{M-\Phi(M)} = \ip{M}{\Phi(M)-\Phi^2(M)} = \ip(M)(0) = 0
\]
and therefore
\[
	\norm{M}^2 = \norm{M-\Phi(M)}^2 + \norm{\Phi(M)}^2.
\]
This implies that $\norm{\Phi(M)}\le\norm{M}$ for any $M$. Hence the operator norm of $\Phi$
is at most 1.

\begin{lemma}\label{lem:PhiDt}
	For any density matrix $D$ and for any time $t$, we have $\Phi(D(t))=\Phi(D)$.
\end{lemma}

\proof
One line:
\[
	\Phi(D(t)) = \sum_r E_rU(t)DU(-t)E_r = \sum_r e^{it\th_r}E_rDE_r e^{-it\th_r}  =\Phi(D).\qed
\]

The \textsl{average mixing matrix} $\hmx$ of the graph $X$ is
\[
	\hmx = \sum_r (E_r)^2.
\]

Our next result is Theorem~3.1 in \cite{us-avmix2017}.

\begin{theorem}\label{thm:mhat-gram}
	If $a,b\in V(X)$, then
	\[
		(\hmx)_{a,b} = \ip{\Phi(D_a)}{\Phi(D_b)}.\qed
	\]
\end{theorem}

If $a\in V(X)$, then
\[
	\Phi(D_a) = \sum_r E_re_ae_a^TE_r
\]
We calculate that 
\[
	\norm{E_re_ae_a^TE_r} = e_a^TE_re_a = (E_r){a,a}
\]
and define
\[
	F_r = \frac1{(E_r)_{a,a}} E_re_ae_a^TE_r.
\]
Thus $F_r$ represents orthogonal projection onto the space of $E_re_a$ and the scalars
\[
	(E_r)_{a,a},\qquad r=1,\ldots,m
\]
are the eigenvalues of $\Phi(D_a)$.

\begin{lemma}
	Assume $a$ and $b$ are vertices in the graph $X$. Then:
	\begin{enumerate}[(a)]
		\item 
		$a$ and $b$ are cospectral if and only the average states $\Phi(D_a)$ and $\Phi(D_b)$ 
		are similar.
		\item
		$a$ and $b$ are strongly cospectral if and only if $\Phi(D_a)=\Phi(D_b)$.
	\end{enumerate}
\end{lemma}

\proof
From Equation~\eqref{eq:phiXvphiX}, we see that $a$ and $b$ are cospectral if and only
if $\Phi(D_a)$ and $\Phi(D_b)$ are. For the second claim we note that $a$ and $b$
are parallel if and only if $E_re_ae_a^TE_r=E_re_be_b^TE_r$ for all $r$, that is, if and only if
the projections $F_r$ are the same for $a$ and $b$.\qed

The sum $\sum_r F_r$ is the matrix representing orthogonal projection onto the walk module
generated by $e_a$.

We introduce spectral densities of subsets of vertices of a graph. 
Assume $S\sbs V(X)$ and let $z$ be the normalized characteristic vector of $S$. 
(So $z$ is zero off $S$, constant on $S$ and $z^Tz=1$.) The quantities
\[
	z^TE_r z,\quad (r=1,\ldots,m)
\]
are non-negative and sum to 1. Hence they determine a probability density on the
eigenvalues of $A$; this is the \textsl{spectral density} of $S$. We will only
work with the case where $S$ is a single vertex, where the value of the spectral density
of vertex $a$ on $\th_r$ is $(E_r)_{a,a}$. Hence the spectral density is determined
by the eigenvalues of $\Phi(D_a)$. The generating function for closed walks on $a$
is the moment generating function for the spectral density at $a$ and, viewed
as a generating function, $U(t)_{a,a}$ is the characteristic function of the
spectral density.

If $\seq p1n$ and $\seq q1n$ are two probability densities on the same finite set,
we define their \textsl{fidelity} to be
\[
	\sum_{j=1}^n \sqrt{p_j}\sqrt{q_j}.
\]
By Cauchy-Schwarz, this is at most 1, and equality holds if and only if $p_j=q_j$
for all $j$. Thus we may view the fidelity as a measure of distance between
probability densities with the same finite support. 

More background on average mixing appears in \cite{cg-average,us-avmix2017}

\section{An Uncomplicated Algebra}

We need information about the matrix algebra generated by $A$ and $e_ae_a^T$ for a vertex $a$.
It is no harder to work with an arbitrary non-zero vector $z$ in place of a vector $e_a$,
so we do.

We use $\alg{S}$ to denote the algebra generated by a set of matrices. The algebra of interest to
is $\alg{A,zz^T}$, where $A$ is an adjacency matrix and $z\in\re^N$.

\begin{lemma}
	Assume $\cA=\alg{A,zz^T}$ for an adjacency matrix $A$ with spectral decomposition
	$A=\sum_r\th_r E_r$. Let $S$ be the set of eigenvalues $\th_r$ such that $E_rz\ne0$.
	If $r\in S$, define
	\[
		F_r = \frac1{z^TE_rz} E_rzz^TE_r,\quad E'_r = E_r-F_r;
	\] 
	if $r\notin S$ then $E'_r=E_r$. Then the matrices
	\[
		E_rzz^TE_s,\ (r,s\in S),
	\]
	together with the non-zero matrices $E'_r$, form a trace-orthogonal basis for $\cA$.
\end{lemma}

\proof
Easy calculations show that the matrices $F_r$ are idempotents ($F_r$ represents orthogonal projection
onto the span of $E_rz$) and they commute with the spectral idempotents. Further $E_kF_r=0$ if $k\ne r$
and $E_rF_r=0$ if $r\notin S$ and $E_rF_r=F_r$ if $r\in S$. One consequence of this is that
the matrices $E'_r$ are pairwise orthogonal and are orthogonal to each matrix $F_s$.

It is also easy to check that distinct matrices of the form $E_rzz^TE_s$ are trace-orthogonal.

Thus it only remains to verify that the given matrices span $\cA$. The key is that
\[
	(A^kzz^TA^\ell)(A^mzz^TA^n) = z^TA^{\ell+m}z\, A^kzz^TA^n
\]
from which it ensues that $\cA$ is spanned matrices of the form $A^kzz^TA^\ell$, along 
with the powers of $A$. The span of the first set of matrices is equal to the span of
the matrices $E_rzz^TE_s$ and the spectral idempotents span the space of polynomials in $A$;
therefore we have an orthogonal basis as claimed.\qed

\begin{corollary}
	If the vertices $a$ and $b$ in $X$ are parallel with the same eigenvalue support,
	then $\alg{A,e_ae_a^T}=\alg{A,e_be_b^T}$.
\end{corollary}

\proof
Suppose $a$ and $b$ are parallel. If $\th_r$ and $\th_s$ lie in the eigenvalue support of $a$ and $b$,
then $E_re_ae_a^TE_s$ and $E_re_be_b^TE_s$ are non-zero scalar multiples of each other. From the
previous lemma it follows that our two algebras are equal.\qed

\begin{corollary}
	Let $X$ be a graph on $n$ vertices and let $a$ and $b$ be parallel vertices in $X$
	with the same eigenvalue support. If the matrix $Q$ commutes with $A$ and $Qe_a=e_a$,
	then $Qe_b=e_b$.\qed
\end{corollary}

\begin{corollary}
	If $a$ and $b$ are strongly cospectral vertices in $X$, then any automorphism of $X$
	that fixes $a$ also fixes $b$.
\end{corollary}

Given this corollary, it is an easy exercise to show that no two vertices in the Petersen
graph are strongly cospectral, but more is true.

The \textsl{characteristic matrix} of a partition $\pi$ is the matrix whose columns are
the characteristic vectors of the cells of $\pi$. If $P$ is the characteristic matrix
of $\pi$, then $P\one=\one$ and $P^TP$ is diagonal with positive diagonal entries.
If $D=(P^TP)^{1/2}$ then the columns of $PD^{-1}$ are pairwise orthogonal unit vectors,
and we call this matrix the \textsl{normalized characteristic matrix} of $\pi$.
We recall that a partition $\pi$ of $V(X)$ is equitable if the column space of $P$ is
$A$-invariant. Alternatively, $\pi$ is equitable if and only if $PD^{-1}P^T$ commutes
with $A$. (Note that $PD^{-1}P^T$ represents orthogonal projection onto $\col(P)$.)

If $\pi$ is a partition of $V(X)$ and $v\in V(X)$, then $\{v\}$ is a cell of $\pi$
if and only if $PD^{-1}P^Te_v=e_v$.

\begin{corollary}
	If $a$ and $b$ are strongly cospectral vertices in $X$ and $\{a\}$ is a cell in
	the equitable partition $\pi$, then $\{b\}$ is also a cell in $\pi$.\qed
\end{corollary}

If $X$ is a graph and $a\in V(X)$, the cells of the distance partition relative
to $a$ are the sets of vertices at a given distance from $a$. It is easy to verify that
if $X$ is strongly regular, then the distance
partition relative to any vertex is equitable. We conclude that if $X$ is strongly
regular and not complete multipartite, no two distinct vertices $X$ in are strongly cospectral.

\section{Eigenspaces and Parallel Vertices}

Our next result provides one way of deciding whether two vertices are parallel.

\begin{lemma}\label{lem:eraerb}
	The projections of $e_a$ and $e_b$ onto the $\theta_r$-eigenspace are parallel
	if and only if $(E_r)_{a,a}(E_r)_{b,b}-(E_r)_{a,b}^2=0$.
\end{lemma}

\proof
Observe that
\[
	(E_r)_{a,b} = e_a^TE_r^TE_re_b = \ip{E_re_a}{E_re_b}
\]
and for any vertex $c$
\[
	(E_r)_{c,c} = \ip{E_re_c}{E_re_c}
\]
whence Cauchy-Schwarz implies that
\[
	(E_r)_{a,b}^2 \le (E_r)_{a,a} (E_r)_{b,b}
\]
with equality if and only if the vectors $E_re_a$ and $E_re_b$ are parallel.\qed

We point out that $(E_r)_{a,a}(E_r)_{b,b}-(E_r)_{a,b}^2$ is the determinant of the
$2\times2$ submatrix of $E_r$ with rows and columns indexed by $a$ and $b$.

If $u$ and $v$ are vertices in $X$, we say an element $f$ in $\re^{V(X)}$ is 
\textsl{balanced} if $f(u)=f(v)$ and is \textsl{skew}\index{skew vector} 
if $f(u)=-f(v)$. A subspace is balanced or skew if each vector in it is balanced or, 
respectively, skew.

\begin{lemma}
	\label{lem:balorskew}
	Two vertices $u$ and $v$ in $X$ are strongly cospectral if and only if each eigenspace
	is balanced or skew relative to the vertices $u$ and $v$.
\end{lemma}

\proof
If $u$ and $v$ are strongly cospectral, then either $E_r(e_u-e_v)=0$ or
$E_r(e_u+e_v)=0$. Since $\col(E_r)$ is the $\theta_r$-eigenspace, it follows that either
each eigenvector in the $\theta_r$-eigenspace is balanced, or each eigenspace is skew.
The converse follows easily.\qed

\begin{lemma}
	Let $S$ be a subset of $V(X)$ such that any two vertices in $S$ are parallel and
	have the same eigenvalue support, of size $s$. Then $|S| \le s$.
\end{lemma}

\proof
Suppose $a\in S$. Denote the non-zero vectors $E_re_a$ by $\seq x1s$. Then for each
vertex $b$ in $S$, we can write $e_b$ as a linear combination of $\seq x1s$.
Since the vectors $e_b$ for $b$ in $S$ are linearly independent, we must have $|S|\le s$.\qed

\section{Parallel Vertices and a Rational Function}

We need an identity due to Jacobi. A proof is given in \cite[Theorem~4.1.2]{cg-blue}.

\begin{theorem}
	\label{thm:jacobi}
	Let $X$ be a graph. If $D\sbs V(X)$, then
	\[
		\det(((tI-A)^{-1})_{D,D}) = \frac{\phi(X\diff D,t)}{\phi(X,t)}.\qed
	\]
\end{theorem}

\begin{corollary}
	\label{cor:polemult}
	Let $\seq \theta1m$ be the distinct eigenvalues of $X$, with corresponding
	spectral idempotents $\seq E1m$. If $D\sbs V(X)$, the multiplicity of $\theta_r$ 
	as a pole of $\phi(X\diff D,t)/\phi(X,t)$ is equal to $\rk((E_r)_{D,D})$.
\end{corollary}

\proof
We have
\[
	((tI-A)^{-1})_{D,D} = \sum_r \frac1{t-\theta_r} (E_r)_{D,D}.
\]
The right side here is the sum of $F=(t-\theta_r)^{-1}(E_r)_{D,D}$ and a matrix $M$ whose 
entries are rational functions with no pole at $\theta_r$. If $n=|V(X)|$, then $\det(F+M)$ 
is the sum of the determinants of the $2^n$ matrices we get from $M$ by replacing each 
subset of its columns by the corresponding subset of columns of $F$. 
This shows that $\rk((E_r)_{D,D})$
is an upper bound on our multiplicity. If $F_r$ were diagonal, we would have equality.
But there is an invertible real matrix $G$ such that $F=G^TDG$ where $D$ is diagonal, with 
nonzero diagonal entries equal to 1. Hence $G^{-T}((tI-A)^{-1})_{D,D}\,G^{-1}$ has a pole of
order $\rk((E_r)_{D,D})$ at $\theta_r$. This completes the proof.\qed

We note that $(E_r)_{D,D}$ is the Gram matrix of the vectors $E_re_u$, for $u$ in $D$.

\begin{lemma}
	\label{lem:parpoles}
	Distinct vertices $a$ and $b$ of $X$ are parallel if and only all poles of the
	rational function $\phi(X\diff\{a,b\},t)/\phi(X,t)$ are simple.
\end{lemma}

\proof
By Corollary~\ref{cor:polemult}, if $D=\{a,b\}$ then the multiplicity of the pole
at $\theta_r$ in $\phi(X\diff D,t)/\phi(X,t)$ is equal to $\rk((E_r)_{D,D})$.
We have
\[
	|(E_r)_{a,b}|^2 = (e_a^T E_re_b)^2 = \ip{E_re_a}{E_re_b}^2 \le \norm{E_re_a}^2\norm{E_re_b}^2
		= (E_r)_{a,a}(E_r)_{b,b}
\]
whence it follows that $\rk((E_r)_{D,D})=1$ if and only if $a$ and $b$ are parallel.\qed

\begin{corollary}
	\label{cor:strcos}
	Distinct vertices $a$ and $b$ of $X$ are strongly cospectral if and only if they
	are cospectral and all poles of $\phi(X\diff\{a,b\},t)/\phi(X,t)$ are simple.\qed
\end{corollary}

One merit of this result is that it enables to decide if two vertices are parallel using
exact arithmetic.

\section{Constructing Strongly Cospectral Pairs}

We present two constructions of strongly cospectral vertices. 

\begin{theorem}
    Let $Z$ be the graph obtained from vertex-disjoint graphs $X$ and $Y$ by joining
    a vertex $u$ in $X$ to a vertex $v$ in $Y$ by a path $P$ of length at least one.
    If $u$ and $v$ are cospectral in $Z$, they are strongly cospectral.
\end{theorem}  

\proof 
Assume $A=A(Z)$ and let $\phi_{u,v}(Z,t)$ denote the determinant of the $uv$-minor 
of $tI-A$. From the spectral decomposition of $A$, we have
\[
    \frac{\phi_{u,v}(Z,t)}{\phi(Z,t)} = ((tI-a)^{-1})_{u,v} 
		= \sum_r \frac{(E_r)_{u,v}}{t-\theta_r},
\]
showing that the poles of $\phi_{u,v}(Z,t)/\phi(Z,t)$ are simple. 
From \cite[Corollary~2.2]{cg-blue}, we have
\[
    \phi_{u,v}(Z,t) = \sum_P \phi(X\diff P,t)
\]                                          
where the sum is over all paths in $X$ that join $u$ to $v$. 
By construction there is only one path in $Z$ that joins $u$ to $v$, and therefore
\[
    \phi_{u,v}(Z,t) = \phi(X\diff u,t)\phi(Y\diff v,t).
\]
If $Q$ is the path we get from $P$ by deleting its end-vertices.
\[
    \frac{\phi(Z\diff\{u,v\},t)}{\phi(Z,t)} 
        = \phi(Q,t)\frac{\phi(X\diff u,t)\phi(Y\diff v,t)}{\phi(Z,t)}
        = \phi(Q,t)\frac{\phi_{u,v}(Z,t)}{\phi(Z,t)}
\] 
We conclude that the poles of $\phi(Z\diff\{u,v\},t)/\phi(Z,t)$ are all simple and so,
by Lemma~\ref{lem:parpoles}, it follows that $u$ and $v$ are strongly cospectral.\qed

Note that $u$ and $v$ will be cospectral in $Z$ if $X$ and $Y$ are cospectral
and also $X\diff u$ and $X\diff v$ are cospectral. We get interesting examples
by taking two vertex-disjoint copies of Schenk's tree from Figure~\ref{fig:Fig1}
and joining the vertex $u$ in the first copy to vertex $v$ in the second by a path
of positive length. This gives pairs of strongly cospectral vertices that do not lie
in an orbit of the automorphism group of the resulting graph.

Now we consider a rabbit-ear construction. Our first step is an interesting
unpublished observation due to K. Guo, reproduced here with her permission.

\begin{lemma}
	If $a$ is a vertex of degree one in $X$ with neighbour $b$, then $a$ and $b$
	are parallel.
\end{lemma}

\proof
Assume $Y=X\diff a$. Then
\[
	\phi(X,t) = t\phi(Y,t) - \phi(Y\diff b,t)
\]
and
\[
	\frac {\phi(X\diff\{a,b\},t)} {\phi(X,t)} 
		= \frac{\phi(X\diff\{a,b\},t)} {t\phi(X\diff a,t)-\phi(X\diff\{a,b\},t)}
		= \frac{1}{t-\frac{\phi(X\diff\{a,b\},t)}{\phi(X\diff a,t)}}.
\]
By interlacing, the derivative of $\phi(X\diff\{a,b\},t)/\phi(X\diff a,t)$ 
is negative wherever it is defined, and therefore the poles of the above rational function
are simple. Now Lemma~\ref{lem:parpoles} implies that $a$ and $b$ are parallel.\qed

We use $\mult(\th,X)$ to denote the multiplicity of $\th$ are a zero if $\phi(X,t)$.

\begin{lemma}\label{lem:rabbit-ear}
	Let $a$ be a vertex in $X$ and let $Z$ be formed from $X$ by joining two new
	vertices of valency one to $a$. If $\mult(0,X\diff a)\le\mult(0,X)$, then the two new 
	vertices are strongly cospectral in $Z$.
\end{lemma}

\proof
Assume the two new vertices are $b$ and $c$. Since $Z\diff b$ and $Z\diff c$ are isomorphic,
$b$ and $c$ are cospectral. We have
\[
	\phi(Z,t) = t^2\phi(X,t)-2t\phi(X\diff a,t)
\]
and so we are concerned with the multiplicities of the poles of
\[
	\frac{\phi(X,t)}{t(t\phi(X,t)-2\phi(X\diff a,t))}
		 = \frac{1}{t\bigl(t-2\frac{\phi(X\diff a,t)}{\phi(X,t)}\bigr)}
\]
By interlacing the zeros of 
\[
	t-2\frac{\phi(X\diff a,t)}{\phi(X,t)}
\]
are simple and hence Lemma~\ref{lem:parpoles} yields that $b$ and $c$ are parallel if and only
if $0$ is not a zero of this rational function. We see that $0$ is a zero if and only if
the multiplicity of $0$ as an eigenvalue of $X\diff a$ is greater than its multiplicity
as an eigenvalue of $X$.\qed

\section{Walk-Regular Graphs}

A graph is \textsl{walk regular} if all its vertices are cospectral. The conept was introduced
in \cite{cgbdm-wlkreg}. Clearly vertex-transitive graphs are walk regular, as a strongly regular graphs.
An old and well-known result states that a vertex-transitive graph with only simple eigenvalues
is $K_1$ or $K_2$. This has been generalized---a walk regular graph with only
simple eigenvalues is $K_1$ or $K_2$ (see e.g., \cite[Theorem~4.8]{cgbdm-wlkreg}). 
The following result generalizes this in turn.

\begin{lemma}
\label{lem:all-strong}
	If all vertices in $X$ are strongly cospectral, then $X=K_2$.
\end{lemma}

\proof
If all vertices of $X$ are strongly cospectral to $u$, then the $\theta_r$ eigenspace
of $X$ is spanned by $E_re_u$, and therefore all eigenvalues of $X$ are simple. 
Assume $n=|V(X)|$. If
\[
	M = \pmat{E_1e_u& \hdots E_ne_u}
\]
then $D=M^TM$ is diagonal. If $S$ is the matrix of coefficients defined above
\[
	S = D^{-1}M^T
\]
and $SS^T=D^{=1}M^TMD^{-1}=D^{-1}$. But $S$ is a $\pm1$-matrix and therefore $SS^T=nI$.
Hence $S$ is a Hadamard matrix and $n$ must be even.

Consequently
\[
	\frac1n = D_{r,r} = e_u^TE_re_u
\]
and therefore each diagonal entry of $E_r$ is equal to $1/n$. It follows that
$X$ is walk-regular and therefore by \cite[Theorem~4.8]{cgbdm-wlkreg} we deduce that $|V(X)|\le2$.\qed

The four vertices of degree two in the Cartesian product of $P_3$ with $K_2$ are pairwise 
strongly cospectral, so we can have more than a pair of strongly cospectral vertices.
(They are cospectral because they form an orbit under the action of the automorphism group.
To see they are parallel it is easiest to note that the characteristic polynomial has only simple
zeros; you can verify this using your favourite computer algebra package.)

\section{Symmetries}

An \textsl{orthogonal symmetry} of a graph $X$ is an orthogonal matrix that commutes
with $A$. If the eigenvalue $\theta_r$ of $X$ has multiplicity $m_r$ and $O(m)$ denotes
the group of $m\times m$ orthogonal real matrices, then the orthogonal symmetries of
$X$ form a group isomorphic to the direct product of the orthogonal groups $O(m_r)$.
Thus this group is determined entirely by the multiplicities of the eigenvalues of $X$
and, given this, does not promise to be very useful. Nonetheless it does have its
applications. Note that the permutation matrices in it form a group isomorphic
to $\aut(X)$.

If the idempotents in the spectral decomposition of $A$ are $\seq E1m$ and
$\sigma_r^2=1$ for each $r$, then
\[
    S = \sum_r \sigma_r E_r
\]
satisfies $S^2=I$. Since $S=S^T$, we see that $S$ is orthogonal. Since $S$ must be
a polynomial in $A$, it follows that the $2^m$ matrices $S$ form a subgroup
of the orthogonal symmetries of $X$; this subgroup is an elementary abelian 2-group.
Any automorphism of $X$ that lies in this group must lie in the centre of $\aut(X)$.

If $a$ and $b$ are cospectral then $A(X\diff a)$ and $A(X\diff b)$ are similar.
Since these matrices are symmetric, there is an orthogonal matrix $L$ say, such that
$L^TA(X\diff a)L=A(X\diff b)$.

\begin{lemma}
    The vertices $a$ and $b$ in $X$ are cospectral if and only there is an orthogonal symmetry
    $Q$ of $X$ such that $Q^2=I$ and $Qe_a=e_b$.
\end{lemma}

\proof
Let $U(+)$ and $U(-)$ respectively denote the $A$-modules generated by $e_a+e_b$
and $e_a-e_b$. By Theorem~\ref{thm:cospvs}(g), these two modules are orthogonal subspaces of $\re^{V(X)}$.
Let $U(0)$ be the orthogonal complement of $U(+)+U(-)$. There is a unique orthogonal
matrix $Q$ such that $Qx=-x$ of $x\in U(-)$ and $Qx=x$ if $x$ lies in $U(+)$ or $U(0)$.

Then
\[
    2Qe_a = Q((e_a+e_b)+(e_a-e_b)) = e_a+e_b -e_a+e_b = 2e_b,
\]
and so $Qe_a=e_b$.

As 
\[
    QA(e_a+e_b) = A(e_a+e_b) = AQ(e_a+e_b)
\]
and as
\[
    QA^k(e_a-e_b) = -A(e_a-e_b) = AQ(e_a-e_b)
\]
If $x\in U(0)$, then
\[
    QAx = Ax = AQx
\]
and therefore $QA=AQ$. 

Thus we have shown that a symmetry exists as required when $a$
and $b$ are cospectral. The converse is straightforward.\qed

It is interesting to note that if $a,b\in V(X)$ and some automorphism $\ga$ maps $a$
to $b$, it does not necessarily follow that $\ga$ maps $b$ to $a$. In fact a permutation
group $G$ on a set $V$ is said to be \textsl{generously transitive} if each pair of elements
of $V$ is swapped by some element of $G$. A transitive group of order cannot be generously 
transitive. The lemma implies that if $\ga$ maps $a$ to $b$, then some orthogonal matrix 
swaps $a$ and $b$, but this matrix need not be related to any automorphism of $X$.

\begin{theorem}
	\label{thm:sc-symm}
    The vertices $a$ and $b$ in $X$ are strongly cospectral if and only there is an 
    orthogonal symmetry $Q$ of $X$ such that $Q$ is a polynomial in $A$, is
    rational, $Q^2=I$ and $Qe_a=e_b$.
\end{theorem}

\proof
We use exactly the same construction as in the previous theorem and then observe that
it $a$ and $b$ are strongly cospectral, then $U(+)$ and $U(-)$ are both direct
sums of eigenspaces of $A$. This implies that $Q$ is a signed sum of the idempotents
$E_r$, and hence is a polynomial in $A$.

Let $\mathbb{E}$ be the extension of the rationals by the eigenvalues of $X$ and let
$\alpha$ be an automorphism of $\mathbb{E}$. Assume $a$ and $b$ are strongly cospectral.
Then $E_r^\alpha$ is an idempotent in the spectral decomposition of $A$, associated
to the eigenvalue $\theta_r^\alpha$. Therefore $((E_r)_{a,a})^\alpha > 0$ and 
consequently $((E_r)_{a,b})$ and $((E_r)_{a,a})^\alpha$ must have the same sign. It
follows that $Q$ is fixed by all field automorphisms of $\mathbb{E}$ and therefore
it is a rational matrix.

The converse is straightforward.\qed

Suppose $X$ is walk regular and $a$ and $b$ are strongly cospectral. Then $Q_{a,a}=0$
but, since $Q$ is a polynomial in $A$, its diagonal is constant. Therefore $\tr(Q)=0$.
Since $Q^2=I$ its eigenvalues are all $\pm1$; we conclude that $1$ and $-1$ have equal
multiplicity and therefore $|V(X)|$ must be even.

With a little more information, we can sharpen the previous theorem and derive a reformulation
of Coutinho \cite[Lemma~3.1(i)]{gc-specex}. Recall that the \textsl{eccentricity}
of a vertex $u$ in $X$ is the least integer $d$ such that any vertex of $X$ is at distance
at most $d$ from $u$. If the eccentricity of $a$ is $d$, then the supports of the
vectors $(A+I)^je_a$ (for $j=0,\ldots,d$) form  a strictly increasing sequence of subsets 
of $V(X)$. Therefore these vectors are linearly independent and accordingly
$d+1$ is a lower bound on the dimension of the walk module $\alg{e_a}_A$. If equality
holds in the bound, Coutinho defines the vertex $a$ to be \textsl{spectrally extremal}.

\begin{corollary}
	Let $a$ and $b$ be strongly cospectral vertices, and assume $a$ has eccentrity $d$.
	If the size of the eigenvalue support of $a$ is equal to $d+1$, then $b$ is the
	unique vertex at distance $d$ from $a$.
\end{corollary}

\proof
Suppose the eigenvalue support of $A$ has size $s$. We have $Qe_a=e_b$ and $Q=p(A)$,
where we choose $p$ to have the least possible degree. It follows that $deg(p)=s-1$.
Since $s-1$ is the eccentricity of $a$, for each vertices $u$ are distance $s-1$
from $a$, the corresponding entry of $p(A)e_a$ is not zero. Therefore $b$ is the
unique vertex in $X$ at distance $s-1$ from $a$.\qed

It can be shown that each vertex in a distance-regular graph is spectrally extremal.

Recall that \textsl{$r$-th distance graph}
$X_r$ of $X$ is the graph with vertex set $V(X)$, where two vertices are adjacent in $X_r$
if thay are distance $r$ in $X$. (Thus $X_1=X$.) We use $A_r$ to denote adjacency matrix of $X_r$
and we set $A_0=I$. We have $\sum_r A_r=J$. We define $X$ to be \textsl{distance regular} if, for
each $r$, the matrix $A_r$ is a polynomial of degree $r$ in $A_1$. It follows from the definition
that $J$ is a polynomial in $A_1$ and consequently $A_r$ and $J$ commute for each $r$. Therefore
the distance graphs $X_r$ are regular.

If $A$ is the adjacency matrix of distance-regular graph, then $A^k$ is a linear combination
of the matrices $\seq A0d$ (for any non-negative integer $k$). Accordingly the diagonal of
$A^k$ is constant for all $k$, and therefore any two vertices in $X$ are cospectral.

We use our theory to present a short proof of a result of Coutinho et al \cite{drg-pst}.

\begin{theorem}
	Suppose $X$ is a distance-regular graph of diameter $d$, with distance matrices $\seq A0d$.
	If $a$ and $b$ are distinct strongly cospectral vertices in $X$, then $A_d$ is a
	permutation matrix of order two and $A_de_a=e_b$.
\end{theorem}

\proof
Let $Q$ be the matrix provided by Theorem~\ref{thm:sc-symm}. Then $Q$ lies in the Bose-Mesner algebra
of the association scheme $\cA=\{\seq A0d\}$ which contains $X$. Since $Qe_a=e_b$, the $a$-column
of $Q$ has exactly one nonzero entry, $Q_{a,b}$. This implies that $Q$ is equal to one
of the matrices $A_r$, and that $A_r$ is a permutation matrix.\qed

A distance-regular graph is \textsl{primitive} if its distance-graphs $\seq X1d$ are connected,
otherwise it is \textsl{imprimitive}. It is a standard result that if a distance-regular graph
of diameter $d$ is imprimitive, either $X_2$ is not connected (and $X$ is bipartite), or
$X_d$ is not connected (in which case the graphs is said to be \textsl{antipodal}). The $d$-cube
is distance-regular, bipartite and antipodal. The previous theorem implies that a distance-regular
graph which contains a pair of strongly cospectral vertices is imprimitive.

\section{Automorphisms, Equitable Partitions}
\label{sec:auto-eqpart}

Let $\pi$ be a partition of $V(X)$. We say that $\pi$ is an \textsl{equitable partition}
if the space of functions on $V(X)$ that are constant on the cells of $\pi$. (There are less
sophisticated definitions, but this one is best suited to our immediate needs. For more
details see, e.g., \cite[Section~9.3]{cggfr-yellow}.) If $Q$ represents orthogonal projection 
onto the space of functions on $V(X)$ constant on the cells of $\pi$, then $\pi$ is equitable 
if and only $A$ and $Q$ commute.

Suppose that we have an equitable partition $\pi$ in which $\{a\}$ is a singleton cell, and
let $Q$ represent orthogonal projection onto the space of functions constant on the cells of
$\pi$. Then $2Q-I$ is orthogonal and commutes with $A$ and $(2Q-I)e_a=e_a$. Now if
$b$ lies in a cell of $\pi$ with size $k$, then 
\[
	\norm{(Q-I)e_b}^2 = (k-1)\frac{1}{k^2} +\left(\frac{1}{k}-1\right)^2 = 1 - \frac1{k}
\]
and so if $Qe_b\ne e_b$, we have $\norm{2(Q-I)e_b}\ge\sqrt{2}$. Therefore:

\begin{lemma}
	Suppose $a,b\in V(X)$. If $\norm{D_a(t)-D_b}<1/\sqrt{2}$, then any equitable
	partition in which $\{a\}$ is a singleton cell must also have $\{b\}$ as a singleton cell.\qed
\end{lemma}

If $D$ is a pure state, then $D^2=D$ and consequently if $D_1$ and $D_2$ are pure states
\[
	\norm{D_1-D_2}^2 = \tr(D_1-D_2)^2 = \tr(D_1+D_2 -2D_1D_2) = 2-2\ip{D_1}{D_2}
\]
If $D_1=yy^*$ and $D_2=zz^*$, this yields that $\norm{D_1-D_2}^2=2-2(y^*z)^2$.

\begin{lemma}
	Let $a$ and $b$ be vertices of $X$. If there is a time $t$ such that 
	$\norm{D_a(t)-D_b}<1/\sqrt{2}$, then any automorphism of $X$ that fixes $a$ must 
	also fix $b$.
\end{lemma}

\proof
Assume $P$ is an orthogonal matrix that commutes with $A$ and $Pe_a=e_a$. Then
\begin{align*}
	P(D_a(t)-D_b)P^T &= PU(t)e_ae_a^TU(-t)P^T -Pe_be_b^TP^T\cr
		&= U(t)Pe_ae_a^TP^TU(-t)-Pe_be_b^TP^T\cr
		&= U(t)D_aU(-t)-Pe_be_b^TP^T\cr
		&= D_a(t) - PD_bP^T.
\end{align*}
This implies that
\[
	\norm{D_a(t)-D_b} = \norm{D_a(t)-PD_bP^T}
\]
and hence if $\norm{D_a(t)-D_b}=\de$, then by the triangle inequality,
\[
	\norm{D_b-PD_bP^T} \le 2\delta.
\]

Now assume $P$ is a permutation matrix. Then $PD_bP^T=D_c$ for some vertex $c$. If $c=b$, then
$Pe_b=e_b$. If $Pe_b\ne e_b$, then
\[
	\norm{D_b-PD_bP^T}^2 = \norm{D_b-D_c}^2 = 2.
\]
We conclude that if there is a time $t$ such that $\norm{D_a(t)-D(t)}<1/\sqrt{2}$, 
then any automorphism of $X$ that fixes $a$ must also fix $b$.\qed

Suppose that we have an equitable partition $\pi$ in which $\{a\}$ is a singleton cell, and
let $Q$ represent orthogonal projection onto the space of functions constant on the cells of
$\pi$. Then $2Q-I$ is orthogonal and commutes with $A$ and $(2Q-I)e_a=e_a$. Now if
$b$ lies in a cell of $\pi$ with size $k$, then 
\[
	\norm{(Q-I)e_b}^2 = (k-1)\frac{1}{k^2} +\left(\frac{1}{k}-1\right)^2 = 1 - \frac1{k}
\]
and therefore if $Qe_b\ne e_b$, we have $\norm{2(Q-I)e_b}\ge\sqrt{2}$. Therefore:

\begin{lemma}
	Assume $a,b\in V(X)$. If $\norm{D_a(t)-D_b}<1/\sqrt{2}$, then any equitable
	partition in which $\{a\}$ is a singleton cell must also have $\{b\}$ as a singleton cell.\qed
\end{lemma}

\section{Cospectral Vertices}

If $(p_i)$ and $(q_i)$ are the spectral densities of two vertices in the graph $X$, then their
fidelity is at most 1, in which case they are equal (and the vertices are cospectral). 
We derive an upper bound on fidelity of the spectral densities of two non-cospectral vertices.
For this we need more machinery.

If $|V(X)|=n$ and $x\in\re^{V(X)}$, the \textsl{walk matrix} of $X$ relative to $x$ is the
$n\times n$ matrix with columns
\[
	x, Ax,\ldots,A^{n-1}x.
\]
The case of combinatorial interest arise when $x$ is the characteristic vector of a nonempty subset
of $V(X)$; in this paper we are concerned only with the case where $x$ is the characteristic vector
of vertex, that is, $x=e_a$ for some vertex $a$. We will use $M_a$ to denote the walk matrix of $X$
relative to the vertex $a$. Note that
\[
	(M_a^TM_a)_{i,j} = e_a^T A^{i+j-2} e_a;
\]
thus the entries of $M_a^TM_a$ are determined by the numbers of closed walks in $X$ that start
(and finish) at $a$.

Suppose 
\[
	A = \sum_r \th_r E_r
\]
is the spectral decomposition of $A$, thus $\seq\th1m$ are the distinct eigenvalues of $A$
and $E_r$ is the matrix that represents orthogonal projection onto the eigenspace belonging
to $\th_r$. Since the spectral idempotents $E_r$ form a basis for the vector space of real
polynomials in $A$, and since $E_r$ is a polynomial in $A$, it follows that the vectors $E_re_a$
span the column space of $M_a$, more precisely, the non-zero vectors $E_re_a$ form an orthogonal
basis for $\col(M_a)$. The set of eigenvalues $\th_r$ such that $E_re_a\ne0$ is
the \textsl{eigenvalue support} of the vertex $a$. (Hence $\rk(M_a)$ is equal to the
size of the eigenvalue support of $a$.)

\begin{lemma}\label{lem:}
	Assume $a$ and $b$ are distinct vertices in the graph $X$ and set $n=|V(X)|$. 
	Let $A=\sum_r\th_rE_r$ be the spectral decomposition of $X$ and let $F$ be the 
	$m\times n$ matrix with $F_{r\ell}=\theta_r^{\ell-1}$. Then
	\[
		\max_r \{\abs{(E_r)_{a,a}-(E_r)_{b,b}}\} \ge \frac{1}{\tr(FF^T)}.
	\]
\end{lemma}

\proof
Let $N_a$ and $N_b$ respectively
denote the $n\times m$ matrices with columns consisting of the vectors $E_re_a$ and $E_re_b$.
If $M_a$ and $M_b$ are the walk matrices of $a$ and $b$ respectively, then
\[
	M_a = N_aF, \quad M_b = N_bF
\]
and
\begin{equation}
	\label{eq:MaTMa}
	M_a^TM_a - M_b^TM_b  = F^T(N_a^TN_a-N_b^TN_b)F
\end{equation}
The matrices $N_a^TN_a$ and $N_b^TN_b$ are diagonal with
\[
	(N_a^TN_a)_{r,r} = (E_r)_{a,a}, \quad (N_b^TN_b)_{r,r} = (E_r)_{b,b}.
\]
Hence
\begin{equation}
	\label{eq:FTNaT}
	F^T(N_a^TN_a-N_b^TN_b)F = \sum_r ((E_r)_{a,a}-(E_r)_{b,b}) F^Te_re_r^TF.
\end{equation}
Let $\eta$ denote the maximum value over $r$ of $\abs{(E_r)_{a,a}-(E_r)_{b,b}}$. Then by the
triangle inequality
\begin{equation}
	\label{eq:EraaErbb}
	\bigl\lVert{\sum_r ((E_r)_{a,a}-(E_r)_{b,b}) F^Te_re_r^TF}\bigr\rVert
		\le \eta \sum_r \norm{F^Te_re_r^TF}.
\end{equation}
We have
\[
	\norm{F^Te_re_r^TF}^2 = \tr(F^Te_re_r^TF\, F^Te_re_r^TF) = (e_r^TFF^Te_r)^2,
\]
whence $\norm{F^Te_re_r^TF}=(FF^T)_{r,r}$ and therefore the right side in \eqref{eq:EraaErbb}
is equal to $\eta\tr(FF^T)$.

If $a$ and $b$ are not cospectral then $M_a^TM_a\ne M_b^TM_b$ and, since these matrices are integer
matrices, the norm of $M_a^TM_a - M_b^TM_b$ is at least 1. So Equations \eqref{eq:MaTMa}, \eqref{eq:FTNaT}
and $\eqref{eq:EraaErbb}$ imply that
\[
	\frac1{\tr(FF^T)} \le \eta.\qed
\]

Our next lemma provides a lower bound on $\abs{U(t)_{a,b}}$.

\begin{lemma}
	If $a,b\in V(X)$ and $\seq E1m$ are the spectral idempotents of $A$, then
	\[
		|(E_r)_{a,a}-(E_r)_{b,b}| < \sqrt{8}\sqrt{1-\abs{U(t)_{a,b}}}
	\]
\end{lemma}

\proof
We have
\[
	U(t)_{a,b} = \sum_r e^{it\th_r} (E_r)_{a,b}.
\]
By the triangle inequality we have
\[
	\abs{U(t)_{a,b}} \le \sum_r \abs{(E_r)_{a,b}}.
\]
Now
\[
	(E_r)_{a,b} = e_a^TE_re_b = \ip{E_re_a}{E_re_b}
\]
and by Cauchy-Schwarz
\[
	\abs{\ip{E_re_a}{E_re_b}} \le \norm{E_re_a}\norm{E_re_b} 
		= \sqrt{(E_r)_{a,a}}\sqrt{(E_r)_{b,b}}.
\]
We conclude that
\[
	\abs{U(t)_{a,b}} \le \sum_r \sqrt{(E_r)_{a,a}}\sqrt{(E_r)_{b,b}}.
\]
Here the upper bound is the fidelity between the spectral densities at $a$ and $b$,
which we denote by vectors $x$ and $t$ respectively. As
\[
	\ip{x-y}{x-y} = 2 - 2\ip xy \le 2 -2\abs{U(t)_{a,b}}.
\]
Therefore, for any $r$, we have
\[
	\biggl\lvert\sqrt{(E_r)_{a,a}} - \sqrt{(E_r)_{b,b}}\biggl\rvert \le \sqrt{2 -2\abs{U(t)_{a,b}}}.
\]
and since, $(E_r)_{a,a}\le1$ and $(E_r)_{b,b}\le1$, we finally have our upper bound:
\[
	|(E_r)_{a,a}-(E_r)_{b,b}| < \sqrt{8}\sqrt{1-\abs{U(t)_{a,b}}}.\qed
\]

We now show that if the orbits of $D_a$ and $D_b$ are close enough, then $a$ and $b$ are cospectral.

\begin{theorem}
	Let $n=V(X)$ and let $\rho$ be the largest eigenvalue of $A$ and let $a$ and $b$ be
	vertices of $X$. If there is a time $t$ such thst
	\[
		\abs{U(t)_{a,b}} \ge 1 - \frac1{8n^4\rho^4},
	\]
	then $a$ and $b$ are cospectral.
\end{theorem}

\proof
We need an estimate for $\tr(FF^T)$. As $\tr(FF^T)$ is equal to the sum of the entries
of the Schur product $F\circ F$, and as the maximum entry of $F$ is $\rho^{n-1}$,
we see that $\tr(FF^T) \le n^2\rho^n$. Now the result follows from the previous two lemmas.\qed

There is a simple relation between $1-\abs{U(t)_{a,b}}$ and the distance between orbits: 
\[
	\norm{D_b-D_a(t)}^2 = 2 - 2\ip{D_b}{D_a(t)} = 2 - 2\abs{U(t)_{a,b}}^2,
\]

\section{Strongly Cospectral Vertices}
\label{sec:near-scospec}

We prove an analog of the result of the previous section, showing that if the orbits
of $D_a$ and $D_b$ are close enough, then $a$ and $b$ are strongly cospectral.

Two preliminary results are needed; the first is Theorem~9.3 in \cite{cg-average}, the second
is Lemma~3.1 from the same source.

\begin{lemma}
	Two vertices of $X$ are strongly cospectral if and only if the corresponding rows of
	$\hmx$ are equal.\qed
\end{lemma}

\begin{lemma}
	Let $D$ denote the discriminant of the minimal polynomial of the adjacency matrix
	of $X$. Then the entries of $D^2\hmx$ are integers.\qed
\end{lemma}

\begin{lemma}
	Let $a$ and $b$ be vertices in the graph $X$. 
	There is a constant $\eta$ (depending on $X$) such that if for some $t$ we have
	\[
		\norm{D_a(t)-D_b} < \eta,
	\]
	then $a$ and $b$ are strongly cospectral.
\end{lemma}

\proof
Suppose $\norm{D_a(t)-D_b}<\zeta$. Then since the operator norm of $\Phi$ is at most $1$,
we can apply Lemma~\ref{lem:PhiDt} to deduce that
\[
	\norm{\Phi(D_a)-\Phi(D_b)} = \norm{\Phi(D_a(t))-\Phi(D_b)} \le \norm{D_a(t)-D_b}<\zeta.
\]
If $u\in V(X)$, then Cauchy-Schwarz yields
\[
	\abs{\ip{\Phi(D_a)-\Phi(D_b)}{\Phi(D_u)}} \le \norm{\Phi(D_a)-\Phi(D_b)}\,\norm{\Phi(D_u)}.
\]
Since $D_u$ is pure, $\norm{D_u}=1$ whence $\norm{\Phi(D_u)}\le1$. and it follows that the 
right side of this inequality is bounded above by $\zeta$. 

We conclude that the absolute value of an entry of $(e_a-e_b)^T\hmx$ is bounded above
by $\zeta$. On the other hand, if $D$ is the discriminant of the minimal polynomial of
$A$, then $D^2\hmx$ is an integer matrix and, accordingly,
if $a$ and $b$ are not strongly cospectral, some entry of $(e_a-e_b)^T\hmx$
is bounded below by $D^{-2}$.\qed

It would not be too difficult to derive an estimate for $\eta$, it would be substantially
smaller than the distance required to show that the vertices are cospectral.

This lemma implies that if there is pretty good state transfer
from $a$ to $b$, then $a$ and $b$ are strongly cospectral.

\section{Problems}

Is there a tree that contains three vertices, any two of which are strongly cospectral?

We have shown that the distance between orbits of $D_a$ and $D_b$ provides a measure
of `similarity' between the vertices $a$ and $b$. Are there further interesting properties
of vertices related to this distance? We admit that computing this distance, even for specific
graphs, is a difficult task. Are there interesting graphs where this computation is feasible?

Find examples of cospectral vertices $a$ and $b$ for which there is a positive constant $\de$
such that $\abs{U(t)_{a,b}}<1-\de$ for all $t$. Find examples of strongly cospectral vertices
satisfying the same condition.

\bibliographystyle{plain}

\end{document}